\documentclass[12pt]{article}

\usepackage{amsmath}
\usepackage{amssymb}

      \newcommand{\RR}{\mathbb R}
      \newcommand{\CC}{\mathbb C}
      \newcommand{\ZZ}{\mathbb Z}

      \newcommand{\qed}{\hfill $\square$}
      \newcommand{\sm}{\mbox{\raisebox{0.2em}{$\scriptstyle \setminus \,$}}}

      \newtheorem{theorem}{Theorem}
      \newtheorem{lemma}{Lemma}
      \newtheorem{prop}{Proposition}
      \newtheorem{coro}{Corollary}

\begin{document}

\begin{center}
\large \bf Symmetries and reversing symmetries of area-preserving
polynomial
          mappings in generalised standard form
\end{center}

\vspace{1.5cm}
\centerline{\large {\sc John A.\ G.\ Roberts}$^{1,2)}$ and {\sc Michael
Baake}$^{3)}$}

\vspace{1.5cm}
{\small
1): Department of Mathematics, La Trobe University, \\
\hspace*{11.2mm}
          VIC 3086, Australia.} 

\bigskip
{\small
2): School of Mathematics, The University of New South Wales, \\
\hspace*{11.2mm}
          Sydney NSW 2052, Australia; \\
\hspace*{11.2mm}
         email: {\sf jagr@maths.unsw.edu.au}

\bigskip
3): Institut f\"{u}r Mathematik und Informatik, Universit\"{a}t Greifswald, \\
\hspace*{11.2mm}
         Friedrich-Ludwig-Jahn-Str. 15a, 17487 Greifswald, Germany; \\
\hspace*{11.2mm}
         email: {\sf mbaake@uni-greifswald.de}
}

\vspace{2cm}
\begin{abstract}
We determine the symmetries
and reversing symmetries within ${\cal G}$, the group of
real planar polynomial automorphisms, of area-preserving  nonlinear
polynomial maps $L$ in generalised standard form,
$$L : \quad x'=x+p^{}_1(y) \, , \; \; y'=y+p^{}_2(x')\, ,$$ where $p^{}_1$
and $p^{}_2$
are polynomials. We do this by using the amalgamated free
product structure of ${\cal G}$. Our results lead to
normal forms for polynomial maps in generalised standard
form and to a classification of the group structures of
the reversing symmetry groups for such maps.
\end{abstract}

\newpage

\section{Introduction  }

In a series of recent articles \cite{baro,baro3,rowi}, the symmetries and
reversing symmetries of some dynamical systems (automorphisms) have been
investigated systematically by means of algebraic methods.
An automorphism $L$ of some space is said
to have a {\em symmetry\/} if there exists an automorphism $S$ satisfying:
\begin{equation} \label{defsym}
        S \circ L \circ S^{-1} = L \, ,
\end{equation}
and a {\em reversing symmetry\/} if there exists an automorphism $R$
satisfying
\begin{equation} \label{defrevsym}
R \circ L \circ R^{-1} = L^{-1}.
\end{equation}
The set of symmetries is non-empty (it certainly contains all powers
of $L$) and this set is actually a group, the {\em symmetry group}
${\cal S}(L)$. On the other hand, the existence {\em a priori} of
any reversing symmetries for a particular $L$ is unknown. When $L$ has a
reversing symmetry, we call it {\em reversible}, and {\em irreversible}
otherwise.

We denote by ${\cal R}(L)$ the set of all symmetries and reversing
symmetries of $L$. It too is a group, called the {\em reversing symmetry
group}
\cite{lamb} of $L$ (see also \cite{good}). If $L$ has no reversing
symmetries or if
$L$ is an involution (i.e. $L^2 = id$), then ${\cal R}(L)={\cal S}(L)$.
On the other hand, if $R$ is a reversing symmetry then so is e.g.
$S^{i} \circ R$ and $R \circ S^{i}$, $i \in \ZZ$
(in particular, $L^{i} \circ R=R \circ L^{-i}$ are reversing symmetries),
whereas $R^{2i}$ is a symmetry.
More technically, if $L$ has reversing symmetries and is not an involution
itself, ${\cal S}(L)$ is a normal
subgroup of ${\cal R}(L)$ with ${\cal R}(L)/{\cal S}(L) \simeq C_2$.
If, in addition, $L$ has an involutory
reversing symmetry $R$ (i.e. $R^2 = id$) then
${\cal R}(L)\simeq{\cal S}(L) \times_s C_2$ (it follows from
(\ref{defrevsym}) that the reversing symmetries $L^{i} \circ R$ and $R
\circ L^{i}$ are then also involutions)\footnote{We denote the cyclic
group of order
$n$ by $C_n$ and use $\times_s$ to denote a semi-direct product.}.
In many cases of reversible automorphisms (and also in the
analogous continuous-time case of reversible flows), it is in fact found
that all reversing symmetries
$R$ satisfying (\ref{defrevsym}) are involutions. In this
case, the automorphism
$L$ can be written as the composition of two involutions, e.g.\
$L \circ R$ and $R$, or $R$ and $R \circ L$. The references \cite{roqu}
and
\cite{laro} include reviews of the properties and applications of reversible
automorphisms and flows.

In the papers \cite{baro,rowi}, the structure of ${\cal R}(L)$ has
been identified when $L$, and possible symmetries $S$ and reversing
symmetries $R$, belong to the group of linear automorphisms of the
two-dimensional torus (or closely related groups). In this case, the
solutions of (\ref{defsym}) and (\ref{defrevsym}) can be solved exactly
using group theoretic and number theoretic techniques. In \cite{baro3},
extensions of some of these results are made for toral automorphisms in
higher dimensions.

In this paper, we consider equations (\ref{defsym}) and
(\ref{defrevsym}) for real polynomial maps $L$ of the plane in
so-called {\em generalised standard form}:
\begin{equation} \label{gsf}
         L: \quad x'=x+p^{}_1(y)\, , \; \; y'=y+p^{}_2(x')\, ,
\end{equation}
where $p^{}_1$ and $p^{}_2$ are (real) polynomials (and the
new $x$ coordinate $x'$ needs to be calculated first so it can be
substituted into the second equation). The form of
(\ref{gsf}) guarantees
that the generalised standard form of map is always area-preserving (even
with
$p^{}_1$ and $p^{}_2$ replaced with non-polynomial functions). As such,
this form
has often been used to create and study models of area-preserving mappings
\cite{mac,mame,roqu}.  Two other standard forms of area-preserving maps,
known as the McMillan and de Vogelaere forms, can always be transformed
to generalised standard form \cite[Section 3.2]{roqu}. For instance, the
McMillan map
\begin{equation} \label{mill}
        M: \quad x'=y\, , \; \; y'=-x+f(y)\, ,
\end{equation}
under the transformation $x \mapsto x$, $y \mapsto x-y$ becomes
(\ref{gsf}) with $p^{}_1(y)=-y$ and $p^{}_2(x')=2x'-f(x')$.

The
much-studied area-preserving H\'enon map in McMillan form is:
\begin{equation} \label{millh}
        H_1: \quad x'=y\, , \; \; y'=-x+2y\, (\lambda+y)\, ,
\end{equation}
and in generalised
standard form is:
\begin{equation} \label{hen}
        H_2: \quad x'=x-y\, , \; \; y'=y+2x'(1-\lambda-x')\, ,
\end{equation}
where $\lambda$ is a parameter. The H\'enon map (\ref{hen}) is the simplest
nonlinear polynomial example of (\ref{gsf}) and all area-preserving
quadratic polynomial mappings can be reduced to it \cite{heno}.
A McMillan map (\ref{mill}) is always reversible with e.g.\ involutory
reversing symmetries: $R: x'=y, \; y'=x$ and $M \circ R:
x'=x, \; y'=-y+f(x)$ (so that $M$ is the composition of involutions
$M \circ R$ and $R$). Hence $H_1$ is certainly reversible and, since
reversibility is invariant under conjugation, $H_2$ is reversible with
e.g.\ involutory polynomial reversing symmetries
   $R: x'=x, \; y'=-y+2x(1-\lambda -x)$ or
$R \circ H_2: x'=x-y, \;
y'=-y$ (so that $H_2$ is the
composition of involutions $R$ and $R \circ H_2$). The
reversibility of $H_1$ or
$H_2$ has been much exploited in numerical studies because it makes it
easier to find periodic orbits \cite{mac,mame,roqu}.

For $L$ of the form (\ref{gsf}), the inverse map $L^{-1}$ is given by
\begin{equation} \label{gsfi}
        L^{-1}: \quad x'=x-p^{}_1(y')\, , \; \; y'=y-p^{}_2(x)\, .
\end{equation}
Hence, if $p^{}_1$ and $p^{}_2$ are polynomial, then $L$ is an example
of a polynomial automorphism of $\RR^2$ to itself since both $L$
and its inverse are polynomial. The polynomial automorphisms
of $\RR^2$ to itself form a group ${\cal G}$ (often called the
general affine group $GA_2$ or the {\em affine Cremona group}
\cite{essen}). In \cite{frmi}, the term {\em generalised
H\'enon transformation} is used for the map
\begin{equation}
      H_g : \; x'=y\, , \; \; y'= -\delta x + p(y),
\end{equation}
with $p$ a polynomial function satisfying $\deg p \ge 2$
and $\delta \neq 0$ the constant Jacobian determinant of $H_g$.
Moreover, it is shown that any nonlinear element of ${\cal G}$ which
is dynamically non-trivial is conjugate within ${\cal G}$ to a composition
\begin{equation} \label{frminorm}
      H_{g^{}_n} \circ H_{g^{}_{n-1}} \circ \ldots H_{g^{}_{2}} \circ 
H_{g^{}_1},
\end{equation}
$n \ge 1$, of generalised H\'enon transformations with $H_{g_i}$ depending
on the constant $\delta_i$ and polynomial $p_i(y)$.
Some of the dynamical properties of
(\ref{frminorm}) have been studied in \cite{frmi,besm,dume}.

Although the
symmetries and reversing symmetries of (\ref{frminorm}) will be studied
more fully elsewhere, the generalised standard form (\ref{gsf}) is a
particular, and significant, area-preserving case of (\ref{frminorm}) when
$n=2$, further justifying its study (see \cite{baro2} and Remark 5 below
in Section 5). Therefore, in this paper, we restrict possible symmetries
$S$ or reversing symmetries
$R$ of $L$ in (\ref{gsf}) to also be polynomial automorphisms
and, with this restriction, we are able to solve (\ref{defsym}) and
(\ref{defrevsym}) exactly.

Our main results, Theorem 1 and its Corollary, are presented in the
next Section.
The key property that allows us to obtain our results in this paper
is that
the group ${\cal G}$ of polynomial automorphisms has a known
structure. This structure is described in Section 3 and its usefulness in
deciding conjugacy problems like (\ref{defsym}) and (\ref{defrevsym})
is pointed out. In Section 4, we give the details behind the proof of
Theorem \ref{main}. Finally, in Section 5, we discuss the implications
of our results, the possible group structures of ${\cal R}(L)$ for
$L$ a polynomial map of generalised standard form, and some extensions
of the work presented here that will appear elsewhere \cite{baro2}.

\section{Statement of results}

We restrict ourselves to the dynamically interesting case of (\ref{gsf})
when $L$ is nonlinear and consider separately the two cases when
both $p^{}_1$ and $p^{}_2$ are nonlinear and when just one of $p^{}_1$ or
$p^{}_2$ is nonlinear. Our main result is:
\begin{theorem} \label{main}
       Let $L$ be a nonlinear polynomial mapping of generalised standard form
       $(\ref{gsf})$. Then $L$ possesses a non-trivial symmetry $S \neq
       L^k, k \in \ZZ$, or a reversing symmetry $R$
       that is also a polynomial automorphism if and only if $p^{}_i$ 
satisfy the
       conditions of Table {\em \ref{tab1}} $($if $\deg p^{}_i \ge 2$,
       $i=1,2\,)$ or Table {\em \ref{tab2}}
       $($if $\deg p^{}_1 =1$, $\deg p^{}_2 \ge 2\,)$.  In
       each case, the corresponding symmetries or reversing symmetries are as
       indicated. The case when $\deg p^{}_2 =1$, $\deg p^{}_1 \ge 2$ follows
from
       Table {\em \ref{tab2}} with replacements $x \rightarrow y$ and
$p^{}_2
\rightarrow -p^{}_1$.
\end{theorem}

\begin{table}
\centerline{
\begin{tabular}{|l|l|}
\hline \vphantom{\Big\|}
        {\em conditions on $p^{}_1$ and $p^{}_2$ of $L$} &
        {\em symmetry (S) / reversing symmetry (R)}
\\ \hline \hline
S1: \; $\begin{array}{l} \vphantom{\Big\|}
             p^{}_1 \mbox{ odd around } e/2 \\ \vphantom{\Big\|}
             p^{}_2 \mbox{ odd around } c/2
             \end{array}$  & $\begin{array}{l}
								S_1:
x'=-x+c, \; \; y'=-y+e  \end{array}$ ($S_1^2 = id$)
\\ \hline
S2: \; $\begin{array}{l} \vphantom{\Big\|}
             p^{}_1 \mbox{ arbitrary } \\ \vphantom{\Big\|}
             p^{}_2(y) = \frac{1}{a} \, p^{}_1(\frac{y}{a}+e)
             \end{array}$ & $\begin{array}{l}
	S_2:x'=ay-ea,\; \; y'=\frac{x+p^{}_1(y)}{a}+e \end{array}$
($S_2^2=L$)
\\ \hline
S3: \; $\begin{array}{l} \vphantom{\Big\|}
             p^{}_1 \mbox{ odd around } \frac{1}{2}(e-\frac{c}{a}) \\
\vphantom{\Big\|}
             p^{}_2(y)=-\frac{1}{a} \,p^{}_1(\frac{y-c}{a})
             \end{array}$ & $\begin{array}{l}
	S_3:x'=ay+c,\; \; y'=-\frac{x+p^{}_1(y)}{a}+e \end{array}$
($S_3^2=S_1 \circ L$)
\\ \hline \hline
R1: \; $\begin{array}{l} \vphantom{\Big\|}
             p^{}_1 \mbox{ arbitrary } \\ \vphantom{\Big\|}
             p^{}_2 \mbox{ odd around } c/2
             \end{array}$  &
	$\begin{array}{l} R_1: x'=-x-p^{}_1(y)+c, \; \; y'=y \quad
(R_1^2 = id)\\
							L \circ R_1:
x'=-x+c, \;\; y'=y-p_2(x) \end{array}$
\\ \hline
R2: \; $\begin{array}{l} \vphantom{\Big\|}
             p^{}_1 \mbox{ odd around } e/2 \\ \vphantom{\Big\|}
             p^{}_2 \mbox{ arbitrary }
             \end{array}$ &
						$\begin{array}{l} R_2:
x'=x+p^{}_1(y), \; \; y'=-y+e \quad (R_2^2 = id)\\
						L \circ R_2: x'=x, \;\;
y'=-y+e+p_2(x) \end{array}$
\\ \hline
R3: \; $\begin{array}{l} \vphantom{\Big\|}
             p^{}_1 \mbox{ even around } e/2 \\ \vphantom{\Big\|}
             p^{}_2 \mbox{ even around } c/2
             \end{array}$ &
$\begin{array}{l} R_3: x'=-x-p^{}_1(y)+c, \; \; y'=-y+e \quad (R_3^2 =
id)\\
            L \circ R_3: x'=-x+c, \;\; y'=-y+e+p_2(x') \end{array}$
\\ \hline
R4: \; $\begin{array}{l} \vphantom{\Big\|}
             p^{}_1 \mbox{ arbitrary } \\ \vphantom{\Big\|}
             p^{}_2(y) = -\frac{1}{a} \,p^{}_1(\frac{y}{a}+e)
             \end{array}$ &
           $\begin{array}{l}R_4: x'=ay-ea, \; \; y'=\frac{x}{a}+e \quad
(R_4^2 =id)\\
						L \circ R_4:
x'=a(y-e-p_2(x)), \; y'=\frac{x}{a}+e+p_2(x')
           \end{array}$
\\ \hline
R5: \; $\begin{array}{l} \vphantom{\Big\|}
             p^{}_1 \mbox{ odd around } \frac{1}{2}(e-\frac{c}{a}) \\
\vphantom{\Big\|}
             p^{}_2(y) = \frac{1}{a}\,p^{}_1(\frac{y-c}{a})
             \end{array}$ & $\begin{array}{l} R_5:
             x'=ay+c,\; \; y'=-\frac{x}{a}+e \end{array}$ ($R_5^4=id$)
\\ \hline
\end{tabular}
}

\caption{Necessary and sufficient conditions for polynomial map $L$ of the
form (\ref{gsf}) with $\deg p^{}_i \ge 2$ to have a non-trivial symmetry
$S$ or reversing symmetry $R$. The parameters $a$, $c$ and $e$ are
arbitrary. In the case of S3 where $S_3^2=S_1 \circ L$, the $S_1$
here is $x'= -x + (c+ae), \; y'=-y+(e-\frac{c}{a})$, since in this case
$p_1$
is odd around $\frac{1}{2}(e-\frac{c}{a})$ and $p_2$ is odd around
$\frac{1}{2}(c+ae)$.
\label{tab1}}
\end{table}

\begin{table}
\centerline{
\begin{tabular}{|l|l|} \hline \vphantom{\Big\|}
        {\em condition on polynomial $p^{}_2$ of $L$} &
        {\em symmetry (S) / reversing symmetry (R)}
\\ \hline \hline
       S1: \; \vphantom{\Big\|}  $\begin{array}{l} p^{}_2
      \mbox{ odd around } e/2 \\ \end{array} $  &
      $\begin{array}{l} S_1:
      x'=-x+e,\; \; y'=-y-2B/A \end{array}$ \;($S_1^2 = id$)
\\ \hline \hline
      R1: \; $\begin{array}{l} p^{}_2 \mbox{ odd around } e/2 \\
      \end{array}$ &
      $\begin{array}{l} R_1:x'=-x+e,\;
      \; y'=y-p^{}_2(x) \quad \quad(R_1^2 = id)
      \vphantom{\Big\|} \\ \vphantom{\Big\|}
      R_1\circ L: x'=-x-p^{}_1(y)+e,\;\;
      y'=y \end{array}$
\\  \hline
      R2: \; $\begin{array}{l} p^{}_2 \mbox{ arbitrary } \\ \end{array}$ &
      $\begin{array}{l} R_2:x'=x,\; \;
      y'=-y+p^{}_2(x)-2B/A \quad (R_2^2 = id)
      \vphantom{\Big\|} \\ \vphantom{\Big\|}
      R_2\circ L: x'=x+p^{}_1(y) \;\;
      y'=-y-2B/A \end{array}$
\\ \hline
\end{tabular}
}

\caption{Necessary and sufficient conditions for polynomial map $L$ of the
form (\ref{gsf}) with $p^{}_1(y)=A y + B$, $\deg p^{}_2 \ge 2$
to have a non-trivial symmetry $S$ or reversing symmetry $R$. The
parameter $e$ is arbitrary.   \label{tab2}}
\end{table}

We reiterate from the Introduction that because the symmetries and
reversing symmetries form a group, other symmetries and reversing
symmetries follow from combining the generating ones presented in Tables
1 and 2 (the groups generated by them are detailed in Section 5,
cf. Tables 5 and 6). In
particular, in the cases of reversibility R1-R4 of Table 1 and R1-R2 of
Table 2, the reversing symmetry $R_i$ in these Tables is easily checked
to be an involution and we also list the involutory reversing symmetry
$L\circ R_i$ or $R_i \circ L$ for the given maps $L$ (with $L$ being the
composition of involutions
$L\circ R_i$ with $R_i$ in Table 1 or $R_i$ with $R_i \circ L$ in Table
2). In the case R5 of reversibility in Table 1, the reversing symmetry
$R_5$ has order 4.  All possible
reversing
symmetries of $L$ are orientation-reversing, except in cases R3 and
R5 of Table \ref{tab1} where
$R$ is orientation-preserving.

Note that a necessary condition for
symmetry possibilities S2 and S3 and reversibility possibilities R4 and
R5 of Table 1 is $\deg p_1 =\deg p_2$. In the case of S2 and S3, the
non-trivial symmetry $S_2$ ($S_3$) is a square root of $L$
($S_1 \circ L$).
More detailed remarks on the results contained in Tables 1 and 2
will be made in Section 5.

Under an affine transformation
\begin{equation} \label{affine}
        x \mapsto \alpha x+\beta \, , \; \; y \mapsto \gamma y+\delta \, ,
\end{equation}
the mapping $L$ of (\ref{gsf}) transforms to $L'$ which is also of generalised
standard form:
\begin{equation} \label{gsft}
        L': \quad x'=x+ \alpha \, p^{}_1 \Bigl(\frac{y-\delta}{\gamma}\Bigr) \,
, \; \;
                  y'=y+ \gamma \, p^{}_2 
\Bigl(\frac{x'-\beta}{\alpha}\Bigr) \, .
\end{equation}
Such transformations can be used to bring the cases of $L$ of
Tables \ref{tab1} and \ref{tab2}
into more transparent normal forms. Since possession of a symmetry or
reversing symmetry is invariant under conjugation, the symmetries
or reversing symmetries of Tables \ref{tab1} and \ref{tab2} transform under
the same
affine transformation used on $L$.

\begin{coro} \label{mainco}
       A polynomial automorphism $L$ of generalised standard form
       $(\ref{gsf})$ possesses a symmetry $S$ or reversing symmetry $R$ that
       is also a polynomial automorphism if and only if, up to an affine
       transformation $(\ref{affine})$, $p^{}_i$ satisfy
       the conditions of Table {\em \ref{tab3}} $($if $\deg p^{}_i \ge 2\,)$ or
       Table {\em \ref{tab4}} $($if $\deg p^{}_1 =1$, $\deg p^{}_2 \ge 2\,)$. In
       each case, the corresponding symmetries or reversing symmetries
       are then as indicated.   \qed
\end{coro}

\begin{table}
\centerline{
\begin{tabular}{|l|l|}
\hline \vphantom{\Big\|}
       {\em conditions on $p^{}_1$ and $p^{}_2$ of $L$} &
       {\em symmetry (S) / reversing symmetry (R)}
\\ \hline \hline
S1: \; $\begin{array}{l} \vphantom{\Big\|}
             p^{}_1 \mbox{ odd around } 0 \\ \vphantom{\Big\|}
             p^{}_2 \mbox{ odd around } 0
             \end{array}$ & $\begin{array}{l} S_1: x'=-x,\; \; y'=-y
   	\end{array}$ ($S_1^2 = id$)
\\ \hline
S2: \; $\begin{array}{l} \vphantom{\Big\|}
             p^{}_1 \mbox{ arbitrary } \\ \vphantom{\Big\|}
             p^{}_2 = p^{}_1
             \end{array}$ & $\begin{array}{l}
	  S_2: x'=y, \; \; y'=x+p^{}_1(y) \end{array}$
             ($S_2^2=L$)
\\ \hline
S3: \; $\begin{array}{l} \vphantom{\Big\|}
             p^{}_1 \mbox{ odd around } 0 \\ \vphantom{\Big\|}
             p^{}_2 = -p^{}_1
             \end{array}$ & $\begin{array}{l}
       	  S_3: x'=y, \; \; y'=-x-p^{}_1(y) \end{array}$
             ($S_3^2=-L$)
\\ \hline \hline
R1: \; $\begin{array}{l} \vphantom{\Big\|}
             p^{}_1 \mbox{ arbitrary } \\ \vphantom{\Big\|}
             p^{}_2 \mbox{ odd around } 0
             \end{array}$  &
	  $\begin{array}{l} R_1: x'=-x-p^{}_1(y), \; \; y'=y \quad
(R_1^2=id) \\
	  L \circ R_1: x'=-x, \;\; y'=y-p_2(x) \end{array}$
\\ \hline
R2: \; $\begin{array}{l} \vphantom{\Big\|}
             p^{}_1 \mbox{ odd around } 0 \\ \vphantom{\Big\|}
             p^{}_2 \mbox{ arbitrary }
             \end{array}$ &
             $\begin{array}{l} R_2: x'=x+p^{}_1(y), \; \; y'=-y
\quad (R_2^2=id) \\
             L \circ R_2: x'=x, \;\; y'=-y+p_2(x) \end{array}$
\\ \hline
R3: \; $\begin{array}{l} \vphantom{\Big\|}
             p^{}_1 \mbox{ even around } 0 \\ \vphantom{\Big\|}
             p^{}_2 \mbox{ even around } 0
             \end{array}$ &
             $\begin{array}{l} R_3: x'=-x-p^{}_1(y), \; \; y'=-y
\quad (R_3^2=id) \\
             L \circ R_3: x'=-x, \;\; y'=-y+p_2(x) \end{array}$
\\ \hline
R4: \; $\begin{array}{l} \vphantom{\Big\|}
             p^{}_1 \mbox{ arbitrary } \\ \vphantom{\Big\|}
             p^{}_2 = -p^{}_1
             \end{array}$ &
             $\begin{array}{l}R_4: x'=y, \; \; y'=x
\quad (R_4^2=id) \\
             L \circ R_4: x'=y+p^{}_1(x), \; y'=x-p^{}_1(x')
             \end{array}$
\\ \hline
R5: \; $\begin{array}{l} \vphantom{\Big\|}
             p^{}_1 \mbox{ odd around } 0 \\ \vphantom{\Big\|}
             p^{}_2 = p^{}_1
             \end{array}$ & $\begin{array}{l} R_5: x'=y,\; \; y'=-x
	  \end{array}$ ($R_5^4=id$)
\\ \hline
\end{tabular}
}

\caption{Normal forms for polynomial map $L$ of the form (\ref{gsf}) with
$\deg p^{}_i \ge 2$ when $L$ possesses a symmetry $S$ or reversing symmetry
$R$ that is a polynomial automorphism. \label{tab3}}
\end{table}

\begin{table}
\centerline{
\begin{tabular}{|l|l|} \hline \vphantom{\Big\|}
        {\em condition on $p^{}_2$ of $L$} &
        {\em symmetry (S) / reversing symmetry (R)}
\\ \hline \hline
S1: \; \vphantom{\Big\|} $\begin{array}{l} p^{}_2 \mbox{ odd around } 0 \\
        \end{array}$ & $\begin{array}{l}
        S_1: x'=-x,\; \; y'=-y \end{array}$ \quad ($S_1^2=id$)
\\ \hline \hline
R1: \; $\begin{array}{l} p^{}_2 \mbox{ odd around } 0 \\ \end{array}$ &
        $\begin{array}{l} R_1:x'=-x,\;
        \; y'=y-p^{}_2(x) \quad (R_1^2 = id)
        \vphantom{\Big\|} \\  \vphantom{\Big\|}
        R_1\circ L: x'=-x-y, \;\; y'=y\end{array}$
\\  \hline
R2: \; $\begin{array}{l} p^{}_2 \mbox{ arbitrary } \\ \end{array}$ &
        $\begin{array}{l} R_2:x'=x,\; \;
        y'=-y+p^{}_2(x) \quad (R_2^2=id)
        \vphantom{\Big\|} \\  \vphantom{\Big\|}
        R_2\circ L: x'=x+y, \; \; y'=-y\end{array}$
\\ \hline
\end{tabular}
}

\caption{Normal forms for polynomial map $L$ of the form (\ref{gsf}) with
$\deg p^{}_1=1$, $\deg p^{}_2 \ge 2$ when $L$ possesses a symmetry $S$
or reversing symmetry $R$ that is a polynomial automorphism.
In each case of the Table, $L$ is given by $x'=x+y,\; \; y'=y+p^{}_2(x')$.
\label{tab4}}
\end{table}

Explicit affine transformations (\ref{affine}) that
convert the cases of Table \ref{tab1} into the corresponding ones of Table
\ref{tab3} are: S1 and R3 ($\alpha=\gamma=1$,
$\beta=- \frac{c}{2} $, $\delta=-\frac{e}{2}$);  S2 ($\alpha=\frac{1}{a}$,
$\beta=e$,
$\gamma=1$, $\delta=0$); S3 and R5 ($\alpha=1$,
$\beta=-\frac{1}{2}(c+ae)$, $\gamma=a$,
$\delta=-\frac{a}{2}(e-\frac{c}{a})$); R1 ($\alpha=\gamma=1$,
$\beta=-\frac{c}{2}$, $\delta=0$); R2 ($\alpha=\gamma=1$, $\beta=0$,
$\delta=-\frac{e}{2}$); R4 ($\alpha=\frac{1}{a}$, $\beta=e$, 
$\gamma=1$, $\delta=0$).
The conversion between Table \ref{tab2} and Table \ref{tab4} is
achieved by (\ref{affine}) with: S1 and R1 ($\alpha=1$,
$\beta=-\frac{e}{2}$, $\gamma=A$, $\delta=B$); R2 ($\alpha=1$, 
$\beta=0$, $\gamma=A$,
$\delta=B$). In all cases, the $p_1$ and $p_2$ of Table \ref{tab3}
and Table \ref{tab4} equal the transformed functions shown in
(\ref{gsft}).

The sufficiency of some of the conditions
in Table \ref{tab3} for allowing a symmetry or reversing symmetry of
(\ref{gsf}), {\em irrespective} of $p_1$ and $p_2$ being polynomials, has been
noted previously in
\cite[Sections 1.2.4.7, 1.1.4.3]{mac} (cases S1, R1 and R2) and \cite{roca}
(case R4). In fact, it is easy to check the sufficiency of all the (reversing)
symmetry conditions of Tables \ref{tab1} -- \ref{tab4}, {\em irrespective} of
$p_1$ and $p_2$ being polynomials.

\section{Algebraic structure of $\mathcal G$}

Let $g$ be a polynomial mapping of $\RR^2$ (or $\CC^2$) of the form:
\begin{equation} \label{genpoly}
        g: \quad x'=P(x,y)\, , \; \; y'=Q(x,y)
\end{equation}
with $P(x,y)$ and $Q(x,y)$ polynomial in $x$ and $y$. Suppose
$g$ is a bijection (i.e.\ one-to-one and onto). Then $g$ has
a well-defined inverse $g^{-1}$. Suppose $g^{-1}$ is also a polynomial
mapping. Then $g$ is called a {\em polynomial automorphism\/} of the real
(or complex) plane.
The set of polynomial automorphisms of the plane forms a group ${\cal G}$ with
a known structure, and this is the key ingredient that we exploit here.
The mathematical context is the theory of amalgamated free products,
convenient sources for background material are \cite[Ch.\ 1.4]{Cohen}
and \cite[Sec.\ 4.2]{mks}

Specifically, let $\cal A$ denote the set of {\em affine} planar maps with
$a \in {\cal A}$ of the form:
\begin{equation} \label{affi}
       a: \quad x'=a^{}_{11}x + a^{}_{12}y + b^{}_1 \, , \; \;
          y'=a^{}_{21}x + a^{}_{22}y + b^{}_2 \, ,
\end{equation}
together with the invertibility condition of its linear part,
\begin{equation}
            a^{}_{11} a^{}_{22} - a^{}_{12} a^{}_{21} \neq 0 \, .
\end{equation}
It is not hard to see that $\cal A$ is a (six dimensional) subgroup of
$\cal G$,
as we have ${\cal A} \simeq K^2 \times_s {\rm GL}(2,K)$ with $K=\RR$ or
$K=\CC$.

Another (infinite-dimensional) subgroup $\cal E$ of $\cal G$ comprises the
so-called {\em elementary} maps
\begin{equation} \label{elem}
         e: \quad x'=\alpha x + p(y) \, , \; \; y'=\beta y + \gamma \, ,
\end{equation}
with $\alpha \beta \neq 0$ and $p(y)$ a polynomial in $y$. Each
$e \in {\cal E}$ has an inverse of similar form (explicitly, $e^{-1}$ is given
by (\ref{elem}) with $\beta$ replaced by $\beta^{-1}$ and $\gamma$ by
$-\beta^{-1} \gamma$).

Notice that the subgroups $\cal A$ and $\cal E$ intersect in those
elementary maps
that are also affine, which occurs whenever $e$ from (\ref{elem}) has
affine $p$,
i.e.\ $p(y) = c^{}_1 y + c^{}_0$. We define the set $\cal I$ by
\begin{equation}
         {\cal I} = {\cal A} \cap {\cal E}
\end{equation}
and note that this intersection is again a subgroup of $\cal G$.

Now, we need the following result of Jung \cite{Jung}, see also \cite[p.\
68]{frmi}.

\begin{theorem} \label{amalg} {\rm [Jung]}
       The group $\cal G$ of polynomial automorphisms of $\,\RR^2$ $($or
$\,\CC^2\,)$
       is the {\em amalgamated free product} of its subgroups $\cal A$ and
$\cal E$,
       understood as groups with real $($complex$\,)$ parameters.
       More explicitly, $\cal G$ is the free product of $\cal A$ and $\cal E$
       amalgamated along their intersection ${\cal I}$, i.e.\
       ${\cal G} = {\cal A} \underset{\cal I}{*} {\cal E}$.   \qed
\end{theorem}

{}For the present purpose, the practical meaning of Theorem \ref{amalg} is
that any $g \in {\cal G}\sm {\cal I}$ can be written as a so-called
{\em reduced word}
\begin{equation} \label{reduce}
        g = g^{}_n \circ g^{}_{n-1} \circ \ldots \circ g^{}_1 \, ,
\end{equation}
where $n\ge 1$ and the $g^{}_i$ alternate in belonging to ${\cal A}\sm
{\cal I}$
or ${\cal E}\sm{\cal I}$. The value of $n \ge 1$ is called the {\em length}
of the reduced word. Note that no reduced word is equal to the identity of
${\cal G}$, see \cite[Cor.\ 2.2]{frmi}. Observe that $g\in{\cal A}\sm{\cal I}$
implies that also $s\circ g$ and $g\circ s$ are in ${\cal A}\sm{\cal I}$, for
arbitrary $s\in{\cal I}$, and analogously for $g\in{\cal E}\sm{\cal I}$.
The algebraic property of amalgamated free products then gives us the
following
structure for reduced words, see \cite[Cor.\ 2.3]{frmi} or
\cite[Thms.\ 25 and 26]{Cohen}.

\begin{prop} \label{normform}
       The word expression $(\ref{reduce})$ for $g\in {\cal G}\sm {\cal I}$ is
       unique up to the modification that we can replace $g^{}_i$, for
any $i > 1$,
       by $g^{}_i \circ s$, with any $s\in {\cal I}$,
       and simultaneously replace $g^{}_{i-1}$ by $s^{-1} \circ g^{}_{i-1}$.
       Consequently, if $g \in {\cal G}\sm {\cal I}$ can be represented both
       by $g^{}_n \circ g^{}_{n-1} \circ \ldots \circ g^{}_1$
       and $g'_{n'} \circ g'_{n'-1} \circ \ldots \circ g'_1$, then $n=n'$ and
       $g^{}_i$ and $g'_i$ both belong to either ${\cal A}\sm {\cal I}$ or
       ${\cal E} \sm {\cal I}$. \qed
\end{prop}

As an example, consider the generalised standard form (\ref{gsf}).
Let the mapping $t$ be given by
\begin{equation} \label{deft}
            t: \quad x'=y \, , \; \; y'=x \, ,
\end{equation}
and, for $i=1,2$, define the following maps:
\begin{eqnarray}
        e^{}_i & : & \quad x'=x+p^{}_i(y) \, , \; \; y'=y \, ,  \nonumber \\
        q^{}_i & : & \quad x'=y \, ,\; \; y'=x+p^{}_i(y) \, , \label{defe} \\
        r^{}_i & : & \quad x'=x \, , \; \; y'=y+p^{}_i(x) \, . \nonumber
\end{eqnarray}
Notice that $q^{}_i=te^{}_i$ (we suppress the symbol $\circ$ from now on).
We see the following:
$t\in {\cal A}\sm {\cal I}$; $e^{}_i \in {\cal E}\sm {\cal I}$ if
$\deg p^{}_i \ge 2$ and $e^{}_i \in {\cal I}$ if $\deg p^{}_i = 1$ or $0$;
$q^{}_i \in {\cal A}\sm {\cal I}$ only if $\deg p^{}_i = 1$ or $0$;
and $r^{}_i \in {\cal A}\sm {\cal I}$ only if $\deg p^{}_i =1$ and $r^{}_i
\in {\cal I}$
if $\deg p^{}_i=0$.

Possible forms of reduced words (\ref{reduce})
for (\ref{gsf}) depend on the nature of $\deg p^{}_i$ and we classify
them as follows: \\
{\bf Type I}: If $\deg p^{}_i \ge 2$, then
\begin{equation} \label{type1}
              L=te^{}_2te^{}_1.
\end{equation}
{\bf Type II}: If $\deg p^{}_2 \ge 2$ and $\deg p^{}_1 \le 1$, then
\begin{equation} \label{type2}
              L=te^{}_2q^{}_1.
\end{equation}
{\bf Type III}: If $\deg p^{}_1 \ge 2$ and $\deg p^{}_2=1$, then
\begin{equation} \label{type3}
              L=r^{}_2 e^{}_1.
\end{equation}
{\bf Type IV}: If $\deg p^{}_1 \ge 2$ and $\deg p^{}_2=0$, then
\begin{equation}
              L: x'=x+p^{}_1(y) \; \; y'=y+K,
\end{equation}
        with $K$ a constant, is itself an element of ${\cal E}\sm {\cal I}$.

{}From (\ref{defsym}) and (\ref{defrevsym}), existence of a
symmetry $S$ (reversing symmetry $R$)
implies a conjugacy between
$L$ and itself (its inverse
$L^{-1}$). If we look for possible (reversing) symmetries that
are themselves polynomial with polynomial inverse, then the
conjugacy question is asked within ${\cal G}$. In this respect,
the known group structure of ${\cal G}$ also helps. We say that
the reduced word (\ref{reduce}) is {\em cyclically reduced} if
$n=1$ or if $g^{}_n$ and $g^{}_1$ do not both belong to ${\cal A}\sm {\cal
I}$ or
${\cal E}\sm {\cal I}$.

Let us recall \cite[Thm.\ 4.6(iii)]{mks}, compare also
\cite[Thm.\ 26]{Cohen} and \cite[Prop.\ 27]{Cohen}.

\begin{prop} \label{conj}
       Let $g$ and $h$ be cyclically reduced elements of ${\cal G}$ that
       are conjugate in ${\cal G}$. If
       $h = h^{}_n h^{}_{n-1} \ldots h^{}_1$ with $n\ge 2$ and
       $h^{}_i$, $h^{}_{i+1}$ as well as $h^{}_n$ and $h^{}_1$ belong to
different
       factors, then $g=s \{ h^{}_n h^{}_{n-1} \ldots h^{}_1\}s^{-1}$ where
       $s\in {\cal I}$ and $\{ h^{}_n h^{}_{n-1} \ldots h^{}_1\}$ means
some cyclic
       permutation of the $n$ letters $h^{}_i$. \qed
\end{prop}

Finally, we also need the following result (with $\RR^* = \RR\setminus\{0\}$).

\begin{lemma} \label{poly}
       Let $p(x)$ be a non-constant real polynomial satisfying the identity:
\begin{equation} \label{rel}
       p(x) = \gamma \, p(\alpha x + \beta),
\end{equation}
       for some $\gamma$,$\alpha \in \RR^*$ and $\beta \in \RR$. Then
       we have:
\begin{itemize}
\item[{\rm (i)}]  if $\alpha=1$, then $\gamma=1$ and $\beta=0$, while no
        restrictions apply to $p(x)$ $($trivial case$\,)$;
\item[{\rm (ii)}] if $\alpha$ is not a root of unity
        $($i.e.\ $\alpha \neq \pm 1 )$,
        then $\gamma \alpha^n=1$ for some $n\ge 1$ and
        $p(x)=c(x+\frac{\beta}{\alpha-1})^n$, with arbitrary $\beta\in\RR$
        and $c\in\RR^*$.
\end{itemize}
The remaining case is $\alpha=-1$ which implies $\gamma^2=1$. In particular:
\begin{itemize}
\item[{\rm (iii)}] if $\alpha=-1$ and $\gamma=1$, then $\beta$ is arbitrary
and
          $p(x)$ is even around $\beta/2$,
          \newline i.e.\ $p(x)=p(\beta/2)+\sum_{n=1}^{N}
c^{}_{2n}(x-\beta/2)^{2n}$,
          $c^{}_{2n} \in \RR$, $\deg p=2N$;
\item[{\rm (iv)}] if $\alpha=\gamma=-1$, then $\beta$ is arbitrary and
          $p(x)$ is odd around $\beta/2$,
          \newline i.e.\ $p(x)=\sum_{n=0}^{N} c^{}_{2n+1}(x-\beta/2)^{2n+1}$,
          $c^{}_{2n+1} \in \RR$, $\deg p=2N+1$.
\end{itemize}
\end{lemma}

{\sc Proof:}
Let $\xi$ be a (real or complex) zero of $p(x)$. Let
$A(x)$ be the affine map $A : x \mapsto \alpha x + \beta$. The identity
(\ref{rel}) implies that $A^{n}(\xi)$ is also a zero of $p(x)$, for all
integer $n \ge 1$. Since $p(x)$ has only finitely many zeros, the orbit
$\{A^{n}(\xi) \mid n \in \ZZ^+\}$ must be
eventually periodic. So we are led to study the finite
orbits under iteration of $\cal A$.

Let $\alpha=1$ in $A(x)$. Then
$A^{n}(\xi)= \xi +n\beta$. All orbits of $A^n$ are infinite and aperiodic
unless $\beta=0$, in which case $A(x)$ is the identity. Substituting
$\alpha=1$ and $\beta=0$ in (\ref{rel}) gives $\gamma=1$ since $p(x)$ is
not the zero polynomial. This gives the trivial case (i).

Suppose $\alpha \neq 1$. Then the dynamics of
$A(x)$ is conjugate to that of the linear mapping
$L: z \mapsto \alpha z$ via the coordinate transformation
$z=x+\frac{\beta}{\alpha-1}$. Since $L^n(z)=\alpha^n z$, we see
that, for $\alpha$ real, there are two types of behaviour. The first occurs
if $\alpha \neq -1$, in which case $L$ has
only one periodic orbit, the fixed point $z=0$, whence
$A$ has only one periodic orbit, the fixed point
$\xi=-\frac{\beta}{\alpha-1}$. It follows that if $\alpha \neq \pm 1$,
$p(x)$ has only one repeated zero over $\CC$,
the real zero $\xi=-\frac{\beta}{\alpha-1}$, so necessarily $p(x)$
takes the form $c(x+\frac{\beta}{\alpha-1})^n$,
where $n \ge 1$ (since $p$ is non-constant by assumption)
and $c \in \RR^*$ is arbitrary. Substituting this form into (\ref{rel})
reveals that we must have $\gamma \alpha^n =1$ and that no restrictions
apply to $\beta$. This gives case (ii).

If $\alpha=-1$, we obtain the second possibility. Then $A(x)=-x+\beta$ is
an involution, so every point has period 2 except for the fixed point at
$\beta/2$. Eq.~(\ref{rel}) can now be written as
$p(\frac{\beta}{2} + x)= \gamma \, p(\frac{\beta}{2} - x)$.
Evaluating this at $x=0$ we see that either $\gamma=1$ or $p(\beta/2)=0$.
If $\gamma=1$, $p$ must be even around $\beta/2$, but no further
restrictions apply. This gives case (iii).
Finally, let $\gamma\neq 1$, so $\xi=\beta/2$ is a root of $p$.
If we write $p$ as a polynomial in powers of $(x-\beta/2)$ we see,
by comparison of coefficients, that only $\gamma=-1$ remains, and that
$p$ must then be odd around $\beta/2$. This gives case (iv). \qed

\smallskip
{\sc Remark}: Lemma \ref{poly} has a natural complex counterpart. Cases
(i) and (ii) are the same, except that one now has to exclude all roots
of unity in (ii). Next, if $\alpha\neq 1$ is a root of unity, i.e.\
$\alpha^n=1$ for some $n\ge 2$, then $\gamma^m = 1$ for some $m|n$ because
iteration of (\ref{rel}) leads to the equation
$$p\Big(\frac{\beta}{1-\alpha} + x\Big) = \gamma^k \,
       p\Big(\frac{\beta}{1-\alpha} + \alpha^k x\Big)$$
for all $k\ge 1$. Consequently, cases (iii) and (iv) are to be replaced
by appropriate cyclic symmetry requirements around $\frac{\beta}{1-\alpha}\,$.

\section{Proof of main result}

To prove our Theorem \ref{main} of Section 2, we  consider the
above-mentioned Type I (\ref{type1}) and Type II (\ref{type2}) word forms
for $L$ in polynomial generalised standard form and use the conjugacy
result Proposition
\ref{conj}. In this Section, we will write our general element of $\cal
I$ as
\begin{equation} \label{defs}
         s: \quad x'=ax+by+c \, , \; \; y'=dy+e \, ,
\end{equation}
with $ad \neq 0$.

\subsection*{Type I words}

Let us first consider the reversibility of Type I words (\ref{type1}), when
both $p^{}_1$ and $p^{}_2$ are nonlinear. We have
\begin{equation}
      L=te^{}_2te^{}_1 \; \; \Rightarrow \; \; L^{-1}=e_1^{-1}te_2^{-1}t,
\end{equation}
noting that $t$ is an involution so that $t^{-1}=t$. The word
for $L$, and hence that for $L^{-1}$, is cyclically reduced.
Assume that $L$ is reversible, so that $L$ and $L^{-1}$ are conjugate
in $\cal G$. Using Proposition \ref{conj}, we must have:
\begin{equation} \label{pos1}
L=te^{}_2te^{}_1=ste_2^{-1}te_1^{-1}s^{-1}
\end{equation}
or
\begin{equation} \label{pos2}
L=te^{}_2te^{}_1=ste_1^{-1}te_2^{-1}s^{-1}.
\end{equation}
Notice that $te_2^{-1}te_1^{-1}$ and
$te_1^{-1}te_2^{-1}$ are the two cyclic permutations of $L^{-1}$
that begin with an element from ${\cal A}\sm{\cal I}$ (by Proposition
\ref{normform} we
must reject the cyclic
permutations that begin with an element from ${\cal E}\sm{\cal I}$ because
this is preserved on adding an $s$ from the left and disagrees
with the word for $L$ which starts with an element of ${\cal A}\sm{\cal I}$).

If (\ref{pos1}) holds, a reversing symmetry is:
\begin{equation} \label{g1}
          R=s e^{}_1
\end{equation}
since (\ref{pos1}) can be rewritten
\begin{equation} \label{pos1r}
L=te^{}_2te^{}_1=ste_2^{-1}te_1^{-1}s^{-1}=se^{}_1(e_1^{-1}te_2^{-1}t)e_1^{-1}s^
{-1}=R L^
{-1} R^{-1},
\end{equation}
and (\ref{defrevsym}) with $R$ replaced by $R^{-1}$ is an equivalent
definition of a reversing symmetry.
Similarly, if (\ref{pos2}) holds, a reversing symmetry is:
\begin{equation} \label{g2}
           R=st
\end{equation}
since (\ref{pos2}) can be rewritten
\begin{equation} \label{pos2r}
L=te^{}_2te^{}_1=ste_1^{-1}te_2^{-1}s^{-1}=st(e_1^{-1}te_2^{-1}t)ts^{-1}=R
L^{-1}
R^{-1}.
\end{equation}
In the case (\ref{g1}), $R$ is nonlinear; in the case (\ref{g2}), $R$ is
affine.

The leading element from ${\cal A}\sm{\cal I}$ on the left hand side
of (\ref{pos1})  is $t$,
whereas the leading element from $A\sm{\cal I}$ on the right hand side is $st$
(the same holds for (\ref{pos2})). Since the left and right hand sides are
different reduced words for the same element $L$, it follows from Proposition
\ref{normform} that:
\begin{equation} \label{sim}
              t=sts',
\end{equation}
where $s'\in {\cal I}$. Writing out (\ref{sim}) readily shows that the
coefficient of $y$ in
the $x'$ part of $s$ and $s'$ vanishes. Consequently, we can take $b=0$ in
$(\ref{defs})$
for $s$ in both (\ref{pos1}) and (\ref{pos2}).

Consider now the identity (\ref{pos1}) in more detail. Using (\ref{deft}),
(\ref{defe}) and (\ref{defs}) with $b=0$, the $x$ component of the right
hand side of (\ref{pos1}) equals
\begin{equation} \label{xcomp}
          x-a\, p^{}_1\Bigl(\frac{y-e}{d}\Bigr),
\end{equation}
and the $y$ component is
\begin{equation} \label{ycomp}
             y-d\,
p^{}_2\left(\frac{x-c}{a}-p^{}_1\Bigl(\frac{y-e}{d}\Bigr)\right).
\end{equation}
Comparing (\ref{xcomp}) to the $x$ component of $L$ in (\ref{gsf})
shows that $p^{}_1(y)$ satisfies:
\begin{equation} \label{con1}
        p^{}_1(y)=-a \, p^{}_1\Bigl(\frac{y-e}{d}\Bigr).
\end{equation}
Comparing (\ref{ycomp}) to the $y$ component of $L$ in (\ref{gsf})
and using (\ref{con1})
shows that $p^{}_2(y)$ satisfies:
$$ p^{}_2\bigl(x+p^{}_1(y)\bigr)=-d \, p^{}_2 \Bigl(
\frac{x+p^{}_1(y)-c}{a}\Bigr), $$
or, equivalently,
\begin{equation} \label{con2}
        p^{}_2(x)=-d \, p^{}_2 \Bigl(\frac{x-c}{a}\Bigr).
\end{equation}

The analysis of the second type of reversibility for a Type I map,
described by (\ref{pos2}),
follows close to that above. Writing out both sides of (\ref{pos2}), again
with $b=0$ in $s$, and equating now leads to the conditions:
\begin{equation} \label{cond1}
            p^{}_1(y) = -a \, p^{}_2\Bigl(\frac{y-e}{d}\Bigr),
\end{equation}
and
\begin{equation} \label{cond2}
            p^{}_2(y) = -d \, p^{}_1\Bigl(\frac{y-c}{a}\Bigr).
\end{equation}
These show that a necessary condition for this type of reversibility is that
$\deg p^{}_1 = \deg p^{}_2$. This can also be seen algebraically from
(\ref{pos2})
and Proposition \ref{normform} since corresponding elements belonging to
${\cal E}\sm{\cal I}$ on the left and right hand sides of (\ref{pos2})
must have the same degree (and $\deg e_i^{-1}=\deg e^{}_i$).

It remains to study the conditions (\ref{con1}) and (\ref{con2}) and
(\ref{cond1}) and (\ref{cond2}). Conditions (\ref{con1}) and (\ref{con2}) are
both examples of the identity $(\ref{rel})$,
and Lemma \ref{poly} can be used to work through the possible values for
$a$,$c$,$d$ and $e$
in $s$ of (\ref{defs}), together with the conditions on $p^{}_1$ and $p^{}_2$.
For example, taking case (i) of Lemma \ref{poly} for $p^{}_1(y)$ in
(\ref{con1}) imposes
$a=-1, d=1$ and $e=0$, which implies in (\ref{con2}) that
$p^{}_2(x)=-p^{}_2(-x+c)$. Hence $p^{}_2$ is odd around $c/2$ for arbitrary
$c$,
$s: x'=-x+c, y'=y$ from (\ref{defs}), and from (\ref{g1}) we obtain the
reversing symmetry $R_1$ shown in R1 of Table \ref{tab1}. In a similar
way, taking
case (iii) or (iv) of Lemma \ref{poly} for $p^{}_1(y)$ gives, respectively,
cases R3 and R2 of Table \ref{tab1}. Taking case (ii) of Lemma \ref{poly}
for
$p^{}_1(y)$
in (\ref{con1}) does not produce a further example of reversibility.
It implies taking $d \neq \pm 1$,
$ad^{-n} =-1$ and
$e$ arbitrary. Substituting into (\ref{con2}) shows that for $p^{}_2$, only
case (ii) of Lemma \ref{poly} is possible, and only if $a \neq \pm 1$
and $da^{-n'}=-1$. The two conditions $ad^{-n} =-1$ and
$da^{-n'}=-1$ imply $(-1)^{n'}\, d^{n n' - 1}=-1$  which cannot
occur since $d \neq \pm 1$ and $n \ge 2$, $n' \ge 2$ by the assumption
$\deg p^{}_i \ge 2$.

The conditions (\ref{cond1}) and (\ref{cond2}) imply
\begin{equation} \label{con3}
       p^{}_1(y) = a d \, p^{}_1\Bigl(\frac{y-(e+dc)}{ad}\Bigr),
\end{equation}
another case of the identity $(\ref{rel})$ in which it is readily
seen case (i) and case (iv) of Lemma \ref{poly} are the only possibilities.
Taking case (i) for $p^{}_1$ in (\ref{con3}) implies $ad=1$, $e+dc=0$
and $p^{}_1$ arbitrary. Substituting in (\ref{cond2}) gives $p^{}_2$ defined
as in R4 of Table \ref{tab1}. The reversing symmetry $R_4$ in this
instance follows
from (\ref{g2}) with $s: x'=ax-ea, y'=y/a + e$. Similarly, taking case (iv) of
Lemma \ref{poly} for (\ref{con3}) gives R5 of Table \ref{tab1}.

Finally, consider the symmetries of $L$ when it is a Type I word
(\ref{type1}). We now need to investigate the necessary
and sufficient conditions for $L=te^{}_2te^{}_1$ to be conjugate to itself in
${\cal G}$. Using Proposition \ref{conj}, we must either have:
\begin{equation} \label{poss1}
        L=te^{}_2te^{}_1=ste^{}_2te^{}_1s^{-1},
\end{equation}
in which case a symmetry is obviously $S=s$, or
\begin{equation} \label{poss2}
          L=te^{}_2te^{}_1=ste^{}_1te^{}_2s^{-1},
\end{equation}
in which case a symmetry is
\begin{equation}
            S=ste^{}_1 \, .
\end{equation}
The analysis in each case is identical to the corresponding
conditions for reversibility
modulo a sign change (stemming from the fact that $e_i^{-1}$ is obtained
from $e^{}_i$ by
$p^{}_i \to -p^{}_i$). For (\ref{poss1}) to hold, the conditions are
\begin{equation} \label{cons1}
        p^{}_1(y) = a \, p^{}_1\Bigl(\frac{y-e}{d}\Bigr) \, , \quad
        p^{}_2(y) = d \, p^{}_2\Bigl(\frac{y-c}{a}\Big) .
\end{equation}
For (\ref{poss2}) to hold, the conditions are
\begin{equation} \label{cons2}
        p^{}_1(y) = a \, p^{}_2\Bigl(\frac{y-e}{d}\Bigr) \, , \quad
        p^{}_2(y) = d \, p^{}_1\Bigl(\frac{y-c}{a}\Bigr) .
\end{equation}
The analysis of (\ref{cons1}) and (\ref{cons2}) is identical to the
corresponding conditions for reversibility and uses
Lemma \ref{poly}.  Conditions (\ref{cons1}) yield S1 of Table \ref{tab1}
(obtained from taking case (iv) of Lemma \ref{poly} for $p^{}_1$)
and conditions (\ref{cons2}) yield S2 and S3 of Table \ref{tab1}.

\subsection*{Type II words}

We consider the (reversing) symmetries of a Type II word (\ref{type2})
in the case that $\deg p^{}_2 \ge 2$ and $\deg p^{}_1 =1$.  Such a
word is not cyclically reduced since $t$ and $q^{}_1$
both belong to $A\sm{\cal I}$. Consequently, we work with such a word
after conjugation with $t$ (and also call this transformed word $L$
for convenience) :
\begin{equation}
        L=t(te^{}_2q^{}_1)t=e^{}_2 r^{}_1 : \quad
        x'=x+p^{}_2(y') \, , \; \; y'=y+p^{}_1(x) \, ,
\end{equation}
which is now cyclically reduced. Conjugation by $t$ involves interchange
of $x$ and $y$ in the original map.
Using Proposition \ref{conj} and Proposition \ref{normform}, we must have:
\begin{equation} \label{pos11}
          L=e^{}_2 r^{}_1 =se_2^{\pm 1}r_1^{\pm 1}s^{-1},
\end{equation}
with the $+$ sign taken when considering symmetries and the $-$ sign
for reversibility. Clearly, a symmetry is given by the affine
transformation $s$ whereas a reversing symmetry is given by
\begin{equation} \label{revs2}
           R=s e_2^{-1}.
\end{equation}

Writing out (\ref{pos11}) with $s$ given by (\ref{defs}) and taking
$p^{}_1=Ax+B$, $A \neq 0$, for the affine function
$p^{}_1$, shows that $b=0$ in $s$ and:
\begin{equation} \label{eq1}
          p^{}_1(ax+c) = \pm d \, p^{}_1(x)
\end{equation}
together with
\begin{equation} \label{eq2}
          p^{}_2(dy+e) = \pm a \, p^{}_2(y) \, .
\end{equation}
Substituting $p^{}_1(x)=Ax+B$ into (\ref{eq1}) gives information about
the parameters of the mapping $s$:
\begin{equation} \label{ss2}
        a = \pm d \, , \; \;  c=(a-1) B / A \, .
\end{equation}
Consequently, (\ref{eq2}) becomes
\begin{equation}
        p^{}_2(y) = \frac{1}{d} \, p^{}_2(dy+e) \, ,
\end{equation}
which is an example of (\ref{rel}) with $\alpha=d=1/\gamma$
and $\beta=e$. Only case (i) ($d=1$, $e=0$ and $p^{}_2$ arbitrary)
and case (iv) ($d=-1$, $e$ arbitrary and $p^{}_2$ odd around $e/2$)
of Lemma \ref{poly} apply. Case (i) gives R2 of Table \ref{tab2}, using
(\ref{ss2}) and (\ref{revs2}) to write out $R_2$, and similarly, case (iv)
gives S1 and R1 of Table \ref{tab2} (remembering to interchange $x$ and $y$
in order to get the results for $te^{}_2q^{}_1$). This solves the problem of
finding symmetries and reversing symmetries for Type II words with $\deg
p^{}_1=1$.

Note that the remaining case of nonlinear maps (\ref{gsf}) after
Type I and II words is when $p^{}_2$ is affine and $p^{}_1$ is nonlinear.
This is a Type III word (\ref{type3})
$L=r^{}_2 e^{}_1$ whose inverse is $L^{-1}=e_1^{-1} r_2^{-1}$. The
(reversing) symmetries of the latter follow from those of $e^{}_2 r^{}_1$,
just considered,
by the replacements $p^{}_2 \rightarrow -p^{}_1$, $p^{}_1 \rightarrow
-p^{}_2$ and
are precisely the (reversing) symmetries of $L=r^{}_2 e^{}_1$ (because
$L$ and $L^{-1}$, for arbitrary $L$, always have the same symmetries,
and have the same reversing symmetries when these are involutions).
This completes the proof of Theorem \ref{main}.

\section{Remarks on results}

In this paper, we have calculated explicitly the symmetries
and reversing symmetries within ${\cal G}$, the group of
planar polynomial automorphisms, of a real nonlinear
polynomial map $L$ in generalised standard form (\ref{gsf}).
The conditions for, and nature of, possible symmetries $S$ and
reversing symmetries $R$ are given in Tables \ref{tab1} and \ref{tab2}. These
results lead to normal forms for polynomial generalised standard
maps with the corresponding (reversing) symmetries in Tables \ref{tab3}
and \ref{tab4}. Some observations about the results in Tables 1-4 were
already made following Theorem \ref{main} and Corollary \ref{mainco}. We
conclude with the following additional remarks:
\vspace{5mm}

{\bf Remark 1:} Suppose $L$ is an element of ${\cal G}$ with a
$k$th root
$M$ also belonging to ${\cal G}$, so $M^k=L$. Then, obviously,
$M$ is a symmetry of $L$ in ${\cal G}$. It follows from  our
results that when $L$ is nonlinear and of generalised standard form, it
can only have a square root in ${\cal G}$ and only
when $p_1$ and $p_2$ are {\em both} nonlinear (because the only
non-trivial symmetry $S_1$ of Tables 2 or 4 is a linear
involution and $(L\circ S_1)^2=L^2$). Case S2 of
Table \ref{tab3} makes clear the only case of a square root that can
occur at the level of the normal form (the conditions at Case S2
of Table \ref{tab1} are the general necessary and sufficient conditions
for this to occur).  When the
conditions in S3 of Table \ref{tab1} or Table \ref{tab3} are satisfied, S1 is
also satisfied and it is seen that $S_3$ is a square root of $S_1 \circ L$
and a fourth root of $L^2$ (since $S_3$
commutes with $S_1$) but neither of the
polynomial automorphisms $S_1 \circ L$ or $L^2$ are in generalised standard
form.
\vspace{5mm}

{\bf Remark 2:} Whenever $L$ is nonlinear and one of $p^{}_1$ or $p^{}_2$
is affine, it follows from R2 of Table \ref{tab2} that $L$ is reversible in
${\cal G}$ since it always
has an involutory reversing symmetry (e.g.\ the H\'enon map in the
form $H_2$ of (\ref{hen})). When both
$p^{}_1$ and $p^{}_2$ are nonlinear, we need to examine the conditions of
Table \ref{tab1} to see whether $L$ is reversible
but it is obvious that irreversibility is the generic case. In
the special case $\deg p^{}_1=\deg p^{}_2=2$,
$L$ is {\em always} reversible since a quadratic polynomial
is always even about some point and case R3 is satisfied
(i.e. $p(x)=ax^2+bx+c= a(x+\frac{b}{2a})^2+ (c-\frac{b^2}{4a})$ so
$p(x-\frac{b}{2a})=p(-x-\frac{b}{2a})$). It turns out then
that the lowest degree combination of $p^{}_1$ and $p^{}_2$ for which
the conditions R$i$, $i=1, \dots , 5$ of Table \ref{tab1} can be avoided
is $\deg p^{}_1=2$, $\deg p^{}_2=3$ or vice versa. After
\cite{frmi}, define the degree of a polynomial map (\ref{genpoly})
by $\deg g=
\max\,(\deg x',\deg y')=\max\,(\deg P(x,y),\deg Q(x,y))$. It is seen
      that for
a map of the form (\ref{gsf}), $\deg L=
\deg p_1 \cdot \,\deg p_2$. Hence the least degree in which
irreversibility in ${\cal G}$ occurs for area-preserving
maps of generalized standard form is $\deg L = 6$ (in \cite{roca},
using the idea of local reversibility, an example is given of
a map (\ref{gsf}) of degree 9 that has no $C^3$ involutory reversing
symmetry). More
generally, for area-preserving maps $g$ in ${\cal G}$ that are not in the
form (\ref{gsf}), degree four is the minimum degree at which
irreversibility occurs (see Remark 5 below).
\vspace{5mm}

{\bf Remark 3:} Tables \ref{tab3} and \ref{tab4} help to make
transparent the various possible structures of
the reversing symmetry group ${\cal R}(L) \subset {\cal G}$.

Firstly, consider the case that $p^{}_1$ in $L$ is affine. From R2
of Table \ref{tab4}, the reversible $L$ can be
taken to be in the form $x'=x+y, \, y'=y+p^{}_2(x')$. If $p^{}_2$ in the
latter is not odd around some value
$e/2$, then
${\cal S}(L) = \langle L \rangle \simeq C_{\infty}$ since S1 of
Table \ref{tab2} is not satisfied and $L$ has no root in ${\cal G}$
by Remark 1 above.
Because $R_2 : x'=x, \, y'=-y+p^{}_2(x)$  of Table \ref{tab4}
is an involution, ${\cal R}(L)$ is a semidirect product
of ${\cal S}(L)$ and the $C_2$ group generated by $R_2$ \cite{baro}
(recall the dihedral group $D_{\infty} \simeq C_{\infty} \times_s
C_2$).
When
$p^{}_2$ is odd around some value $e/2$, the symmetry group increases
to accommodate the involutory symmetry $S_1$ of Table
\ref{tab2}. Note that
$R_1$ of  Table \ref{tab2} now becomes an additional reversing
symmetry of
$L$ and is the composition of $S_1$ with $R_2$. In fact, in this
case,  $R_2$
and $S_1$ commute and we can simplify the structure of
${\cal R}(L)$ to that of a direct product of its $D_{\infty}$ subgroup
$(\langle L\rangle \times_s \langle R_2\rangle)$ and its $C_2$ subgroup
$\langle S_1\rangle$. The
situation is summarized in Table
\ref{tab5}.

Secondly, consider the case that $p^{}_1$ and $p^{}_2$ are both
nonlinear. Table 6 summarizes the possible distinct symmetry
group and reversing symmetry group
structures and their generators. The top case of Table 6 corresponds to
the absence of reversing symmetries {\em and\/} nontrivial symmetries.
Whenever S2 or S3 of Table \ref{tab1}
hold in Table 6, the $C_\infty$ part of ${\cal S}(L)$ is generated by the
corresponding $S_2$ or $S_3$ of Table \ref{tab1} -- otherwise the
$C_\infty$ part is generated by $L$ itself. Table 6 illustrates various
ways that the cases of Table \ref{tab1} or Table \ref{tab3} interrelate.
For example, it is seen from Table \ref{tab1} that whenever the
conditions for S1 hold, then R1 and R2 follow automatically and
the reversing symmetry $R_2$   is seen to be the composition
of the $R_1$   together with the symmetry $S_1$. Also,
whenever the conditions for R5 hold in Table \ref{tab1}, then
S1, S2, R1 and R2 immediately follow. Furthermore, since a
non-zero polynomial cannot be even about some point and odd
about some (other) point, it follows that reversibility possibility
R3 of Table \ref{tab1} or Table \ref{tab3} can only co-exist with
R4 or S2, with any two of R3, R4 and S2 implying the third.

Recall that a necessary condition for
symmetry possibilities S2 and S3 and
reversibility possibilities R4 and R5 of Table \ref{tab1} is $\deg p_1
=\deg p_2$. It is also clear that symmetry possibilities S2 and S3
cannot co-exist. As found in the S1 possibility of Table \ref{tab5},
it is possible to simplify the structure of ${\cal R}(L)$ in Table
\ref{tab6} to a
direct product for
possibility S1;not S2 and possibility S3. This follows because $R_1$ and
$S_1$ commute in these cases and, significantly in the case S3, because
$R_1$ is a reversing symmetry of $S_3$ (as well as being a reversing
symmetry of $L=S_1 S_3^2$). Note in possibility S3, the reversing
symmetry $R_4=S_3 \circ R_1$. Interestingly, in  possibility S1
and S2 of Table \ref{tab6} we {\em cannot}
achieve ${\cal R}(L)$ as a direct product because
$R_1$ is a reversing symmetry of
$L=S_2^2$ but not of $S_2$ (in fact, $R_1 \circ S_2 \circ R_1 =
S_1 \circ S_2^{-1}$). In
other words, $\langle S_2\rangle$ is not a normal subgroup of
$\langle S_2, R_1\rangle$. The map $S_2$ is
said to have $R_1$ as a reversing 2-symmetry because $R_1$
conjugates $S_2^2$ to its inverse but not $S_2$ to its inverse.
The reversing symmetry $R_5$ in possibility S1 and S2 is obtained
from $R_5=S_2 \circ R_1$.

\begin{table}
\centerline{
\begin{tabular}{|c|c|c|} \hline \vphantom{\Big\|}
        {\em ${\cal S}(L)$} & {\em ${\cal R}(L)$} & {\em Possibilities}
\\ \hline \hline   \vphantom{\Big\|}
        $C_{\infty}$ &  $C_{\infty} \times_s C_2 \simeq D_{\infty}$ &
\\  \vphantom{\Big\|}
    $\langle L\rangle$ &  $\langle L\rangle \times_s \langle R_2\rangle
$ &  not S1
\\ \hline \vphantom{\Big\|}
        $C_{\infty} \times C_2$  &  $(C_{\infty} \times C_2) \times_s C_2
        \simeq D_\infty \times C_2 $ &
\\ \vphantom{\Big\|}
     $\langle L\rangle \times \langle S_1\rangle$  &
     $(\langle L\rangle \times_s \langle R_2\rangle) \times \langle
S_1\rangle$ & S1
\\ \hline
\end{tabular}
}

\caption{Possible group structures of ${\cal S}(L)$ and
${\cal R}(L)$ with their generators when $L$ of
(\ref{gsf}) has one of $p^{}_1$ or $p^{}_2$ affine.  The term `S1'
or `not S1' refers to satisfaction or non-satisfaction of the
condition listed for S1 in Table \ref{tab2}. \label{tab5}}
\end{table}

\begin{table}
\centerline{
\begin{tabular}{|c|c|l|} \hline \vphantom{\Big\|}
        {\em ${\cal S}(L)$} & {\em ${\cal R}(L)$} & \quad {\em Possibilities}
\\ \hline \hline \vphantom{\Big\|}
      $C_{\infty}$ & ${\cal S}(L)$ &
\\ \vphantom{\Big\|}
      $\langle L\rangle$ & $\langle L\rangle$ &
      $\mbox{not } {\rm S}i, \,
        i=1,2,3; \mbox{ not } {\rm R}i, \, i=1, \dots , 5 $
\\ \vphantom{\Big\|}
     $\langle S_2\rangle$ & $\langle S_2\rangle$ &
         $ {\rm S}2; \mbox{ not } {\rm R}i, \, i=1, \dots, 5$
\\ \hline \vphantom{\Big\|}
     $C_{\infty}$ &  $C_{\infty} \times_s C_2 \simeq D_{\infty}$ &
\\ \vphantom{\Big\|}
     $\langle L\rangle$ &  $\langle L\rangle \times_s \langle R_i\rangle $ &
     $ \mbox{not } {\rm S}i, \, i=1,2,3;
        \mbox{ one } {\rm R}i, \, i=1, \dots, 4  $
\\ \vphantom{\Big\|}
     $\langle S_2\rangle$ &  $\langle S_2\rangle \times_s \langle R_3\rangle$ &
     ${\rm S}2 \mbox{ and } {\rm R}3 \;
(\Leftrightarrow {\rm S}2 \mbox{ and } {\rm R}4)$
\\ \hline \vphantom{\Big\|}
     $C_{\infty} \times C_2 $  &
     $\begin{array}{c} \vphantom{\Big\|} (C_{\infty} \times C_2) \times_s C_2 \\
        \simeq D_\infty \times C_2  \end{array}$ &
\\ \vphantom{\Big\|}
     $\langle L\rangle \times \langle S_1\rangle$  &
     ($\langle L\rangle \times_s \langle R_1\rangle)
     \times \langle S_1\rangle $ &
     ${\rm S}1 ; \, \mbox{not } {\rm S}2 \;
     (\Rightarrow {\rm R}1 \mbox{ and } {\rm R}2) $
\\ \vphantom{\Big\|}
     $\langle S_3\rangle \times \langle S_1\rangle$  &
     ($\langle S_3\rangle \times_s \langle R_1\rangle) \times
     \langle S_1\rangle $ &
        $ {\rm S}3 \; (\Rightarrow {\rm S}1, \, {\rm R}1, \, {\rm R}2
        \mbox{ and } {\rm R}4) $
\\ \hline \vphantom{\Big\|}
     $C_{\infty} \times C_2$  &
     $(C_{\infty} \times C_2) \times_s C_2  $ &
\\ \vphantom{\Big\|}
     $\langle S_2\rangle \times \langle S_1\rangle$  &
     ($\langle S_2\rangle \times \langle S_1\rangle) \times_s
       \langle R_1\rangle $ &
        ${\rm S}1 \mbox{ and } {\rm S}2 \;
        (\Leftrightarrow {\rm R}5; \; \Rightarrow {\rm R}1
        \mbox{ and } {\rm R}2) $
\\ \hline
\end{tabular}
}

\caption{Possible group structures of ${\cal S}(L)$ and
${\cal R}(L)$ with their generators when $L$ of (\ref{gsf})
has both $p^{}_1$ and $p^{}_2$ nonlinear. The term `S1' or `not S1' etc.\
refers to the satisfaction or non-satisfaction of the condition listed
for S1 in Table \ref{tab1}. \label{tab6}}
\end{table}

\vspace{5mm}
{\bf Remark 4:}
Veselov \cite{ves1,ves2} has shown that an area-preserving
map $L$ in ${\cal G}$ 
possesses a nonconstant polynomial or rational
integral $I(x,y)=I(x',y')$ if and only if $L$ has
a symmetry $S$ with 
$\langle L, S\rangle 
\simeq C_\infty \times C_\infty$.
>From Tables 5 and 6, this symmetry
group structure does not occur for the polynomial generalised
standard maps considered here so integrable cases are absent.

\vspace{5mm}
{\bf Remark 5:} In future work \cite{baro2}, we will consider ${\cal R}(L)$
for arbitrary $L \in {\cal G}$, i.e. $L$ not necessarily in generalised
standard form. There we can show that some of the results found
here hold more generally. Furthermore, we can obtain these results using
just algebraic properties of the amalgamated free product structure
(rather than by means of constructive calculations for the standard form
of $L$, as we did here).
For instance, if $L \in  {\cal G}$ is not conjugate to an affine
map or to an elementary map, it can be shown
that ${\cal S}(L) \subset {\cal G}$ is always an Abelian group of
the form $C_\infty$ or $C_\infty \times C_2$. Furthermore, any reversing
symmetries $R \in {\cal G}$ are involutions or of order 4.

A necessary condition for $L \in {\cal G}$ to be reversible is that
$L$ is area-preserving or anti area-preserving (i.e. $\det dL=-1$).
To study reversibility, one can use
from \cite{frmi} the
normal form (\ref{frminorm}) for a nonlinear (non-elementary)
$L \in {\cal G}$.
The degree of the normal form (\ref{frminorm}) can be shown \cite{frmi} to
be
$\prod_i \deg p_i$, with $\deg p_i \ge 2$.
When (\ref{frminorm}) has $n=1$ and $\delta_1 = 1$, the area-preserving
mapping $H_{g^{}_1}$ is precisely of McMillan form (\ref{mill}) and
so is reversible. In particular, this implies whenever $\deg L$ is
prime and $L$ is area-preserving, then $L$ is reversible in ${\cal G}$.
If we consider the area-preserving case of (\ref{frminorm}) with
$n=2$ (i.e. $\delta_1 \delta_2=1$), we find that $H_{g^{}_2} \circ
H_{g^{}_1}$ is reversible when $\deg p_1 \neq \deg p_2$ if and only if:
(i) $\delta_1 =\delta_2 = 1$ (in which case $H_{g^{}_i}$ are both
(reversible) McMillan maps of the form (\ref{mill}) and $H_{g^{}_2} \circ
H_{g^{}_1}$ always has e.g. the
involutory reversing symmetry $R: x'=-x+p_1(y),\; y'=y$) ; or
(ii) $\delta_1 =\delta_2 = -1$ (i.e. $H_{g^{}_2} \circ
H_{g^{}_1}$ is in the generalised
standard form (\ref{gsf})) and at least one of the conditions
R1, R2 or R3 of
Table \ref{tab1} is satisfied. In particular, with respect to (ii), if
$\deg p_1$ is odd and $\deg p_2$ is even or vice versa, 
Table
\ref{tab1} shows that there are no non-trivial symmetries and the only
possibility for reversibility is one of R1 and R2.

When $\deg p_1 = \deg p_2$,
an area-preserving $H_{g^{}_2} \circ
H_{g^{}_1}$ has, in addition to (i) and (ii), other
possibilities to be reversible but these are not generically satisfied.
This is true even for the least possible case of $\deg \,(H_{g^{}_2} \circ
H_{g^{}_1})
=4$ (i.e. $\deg p_1 =\deg p_2 = 2$) when the additional reversibility
condition for $p_1(y)=y^2+C_1$, $p_2(y)=y^2+C_2$ in the
$(\delta_1,C_1,C_2)$-space is: (iii) $C_1= \delta_1^{2/3} C_2$
(this form for $p_1$ and $p_2$ in $H_{g^{}_2} \circ
H_{g^{}_1}$ can be taken
without loss of generality from
\cite[Theorem 2.6]{frmi}). Thus, considering  quartic polynomial maps
in ${\cal G}$ that are area-preserving, these are always
reversible if they are derived from (\ref{frminorm}) with $n=1$
or are derived from (\ref{frminorm}) with $n=2$ and satisfy
$\delta_1=\delta_2=\pm 1$ (noting, from Remark 2 above, that
(\ref{gsf}) is always reversible when $\deg L = 4$).
Otherwise, if they are derived from
(\ref{frminorm}) with $n=2$, $\delta_1 \delta_2=1$ and
$\delta_1 \neq \pm 1$, they are generically irreversible in ${\cal G}$.

\section*{Acknowledgements}
This work was supported by a La Trobe University Central Starter Grant
and by the German Research Council (DFG).


\bigskip

\begin{thebibliography}{99}

\bibitem{baro}
M.~Baake and J.\thinspace A.\thinspace G.~Roberts,
Reversing symmetry group of $GL(2,\ZZ)$ and $PGL(2,\ZZ)$
matrices with connections to cat maps and trace maps,
{\em J.\ Phys.} {\bf A 30} (1997) 1549--1573.

\bibitem{baro2}
M.~Baake and J.\thinspace A.\thinspace G.~Roberts,
Reversing symmetry group of polynomial automorphisms of $\RR^2$,
in preparation.

\bibitem{baro3}
M.~Baake and J.\thinspace A.\thinspace G.~Roberts,
Symmetries and reversing symmetries of toral automorphisms,
{\em Nonlinearity} {\bf 14} (2001) R1--R24.

\bibitem{besm}
E.~Bedford and J.~Smillie,
Polynomial diffeomorphisms of $\CC^2$. V. Critical points
and Lyapunov exponents,
{\em J.\ Geom.\ Anal.\ } {\bf 8} (1998) 349--383.

\bibitem{Cohen}
D.\thinspace E.~Cohen,
{\em Combinatorial Group Theory: A Topological Approach},
(Cambridge University Press, Cambridge, 1989).

\bibitem{dume}
H.\thinspace R.~Dullin and J.\thinspace D.~Meiss,
Generalized H\'enon maps: the cubic diffeomorphisms of the
plane,
{\em Physica} {\bf D 143} (2000) 262--289.

\bibitem{essen}
A.~van den Essen,
Seven lectures on polynomial automorphisms,
in: {\em Automorphisms of Affine Spaces}, ed.\ A.\ van den Essen
(Kluwer, Dordrecht, 1995) pp. 3--39.

\bibitem{frmi}
S.~Friedland and J.~Milnor,
Dynamical properties of plane polynomial automorphisms,
{\em Ergod.\ Th.\ \& Dyn.\ Syst.} {\bf 9} (1989) 67--99.

\bibitem{good}
G.~R.~Goodson,
Inverse conjugacies and reversing symmetry groups,
{\em Amer.\ Math.\ Monthly} {\bf 106} (1999) 19--26.

\bibitem{heno}
M.~H\'enon,
Numerical study of quadratic area preserving mappings,
{\em Quart.\ Appl.\ Math.} {\bf 27} (1969) 291--312.

\bibitem{Jung}
H.\thinspace W.\thinspace E.~Jung,
\"Uber ganze birationale Transformationen der Ebene,
{\em J.\ Reine Angew.\ Math.} {\bf 184} (1942) 161--174.

\bibitem{lamb}
J.\thinspace S.\thinspace W.~Lamb,
Reversing symmetries in dynamical systems,
{\em J.\ Phys.} {\bf A 25} (1992) 925--937.

\bibitem{laro}
J.\thinspace S.\thinspace W.~Lamb and J.\thinspace A.\thinspace G.~Roberts,
Time-reversal symmetry in dynamical systems: A survey,
{\em Physica} {\bf D 112} (1998) 1--39.

\bibitem{mac}
R.\thinspace S.~MacKay,
{\em Renormalisation in Area-Preserving Maps}
(World Scientific, Singapore, 1993), and references therein.

\bibitem{mame}
R.\thinspace S.~MacKay and J.\thinspace D.~Meiss (eds.),
{\em Hamiltonian Dynamical Systems}
(Adam Hilger, Bristol, 1987).

\bibitem{mks}
W.~Magnus, A.~Karrass, and D.~Solitar,
{\em Combinatorial Group Theory}, 2nd ed.,
(Dover, New York, 1976).

\bibitem{roca}
J.\thinspace A.\thinspace G.~Roberts  and H.\thinspace W.~Capel,
Area preserving mappings that are not reversible,
{\em Phys.\ Lett.} {\bf A 162} (1992) 243--248.

\bibitem{roqu}
J.\thinspace A.\thinspace G.~Roberts and G.\thinspace R.\thinspace W.~Quispel,
Chaos and time-reversal symmetry --- order and chaos in reversible
dynamical systems,
{\em Phys.\ Rep.} {\bf 216} (1992) 63--177.

\bibitem{rowi}
J.\thinspace A.\thinspace G.~Roberts and R.\thinspace S.~Wilson,
Reversibility of orientation-reversing cat maps and the amalgamated
free product structure of PGL(2,$\ZZ$),
preprint (2002).

\bibitem{ves1}
A. \thinspace P.~Veselov,
Integrable maps, 
{\em Russian Math.\ Surveys} {\bf 46} (1991) 1--51.

\bibitem{ves2}
A. \thinspace P.~Veselov,
Growth and integrability in the dynamics of mappings, 
{\em Commun.\ Math.\ Phys.} {\bf 145} (1992) 181--193.
\end{thebibliography}
\end{document}